\documentclass[twoside, 10pt]{article}
\usepackage{amssymb, amsmath, mathrsfs, amsthm}
\usepackage{graphicx}
\usepackage{mathtools}
\usepackage{color}
\usepackage{float, caption, subcaption}

\newcommand{\bd}{\begin{description}}
\newcommand{\ed}{\end{description}}
\newcommand{\bi}{\begin{itemize}}
\newcommand{\ei}{\end{itemize}}
\newcommand{\be}{\begin{enumerate}}
\newcommand{\ee}{\end{enumerate}}
\newcommand{\beqs}{\begin{eqnarray*}}
\newcommand{\eeqs}{\end{eqnarray*}}

\definecolor{DarkGreen}{rgb}{0.2, 0.6, 0.3}

\newtheorem{theorem}{Theorem}
\newtheorem{lemma}{Lemma}
\newtheorem{definition}{Definition}

\newtheorem{case}{Case}
\newtheorem{subcase}{Subcase}[case]
\newtheorem{claim}{Claim}

\newtheorem{fact}{Fact}

\begin{document}
\title{\bf Gallai-Ramsey number of an 8-cycle}
\author{Jonathan Gregory\footnote{Department of Mathematical Sciences, Georgia Southern University, Statesboro, GA 30460, U.S.A.}, Colton Magnant\footnote{Department of Mathematics, Clayton State University, Morrow, GA 30260, U.S.A}, Zhuojun Magnant\footnote{Department of Mathematical Sciences, Georgia Southern University, Statesboro, GA, 30460, U.S.A.}}
\maketitle

%------------------------abstract------------------------------------------------------
\begin{abstract}
Given graphs $G$ and $H$ and a positive integer $k$, the Gallai-Ramsey number $gr_{k}(G : H)$ is the minimum integer $N$ such that for any integer $n \geq N$, every $k$-edge-coloring of $K_{n}$ contains either a rainbow copy of $G$ or a monochromatic copy of $H$. These numbers have recently been studied for the case when $G = K_{3}$, where still only a few precise numbers are known for all $k$. In this paper, we extend the known precise Gallai-Ramsey numbers to include $H = C_{8}$ for all $k$.
\end{abstract}

%%%%%%%%%%%%%%%%%%%%%%%%%-Section 1-%%%%%%%%%%%%%%%%%%%%%%%%%%%%%%%%%%%%%%%%%%%%%%%%%%%%%%
\section{Introduction}

In this work, we consider only edge-colorings of graphs.  A coloring of a graph is called \emph{rainbow} if no two edges have the same color.  

Colorings of complete graphs which contain no rainbow triangle have very interesting and somewhat surprising structure. In $1967$, Gallai \cite{G67} first examined this structure under the guise of transitive orientations.  This result was restated in \cite{GS04} in the terminology of graphs and can also be traced back to \cite{CE97}.  For the following statement, a \emph{trivial partition} is a partition into only one part.

\begin{theorem}[\cite{G67,GS04}]
\label{Gallai}
In any coloring of a complete graph with at least $2$ vertices containing no rainbow triangle, there exists a non-trivial partition of the vertices (called a \emph{Gallai partition}) such that there are at most two colors on the edges between the parts and only one color on the edges between each pair of parts.
\end{theorem}

In honor of this result, rainbow triangle-free colorings have been called \emph{Gallai colorings}. The partition given by Theorem~\ref{Gallai} is called a \emph{Gallai partition} or \emph{G-partition} for short. Given a Gallai coloring of a complete graph and its associated G-partition, define the \emph{reduced graph} of this partition to be the induced subgraph consisting of exactly one vertex from each part of the partition.  Note that the reduced graph is a $2$-colored complete graph.

When considering $2$-colored complete graphs, a very natural problem to consider is the Ramsey problem of finding a monochromatic (one-colored) copy of some desired subgraph.  Given a graph $G$, let $R_{k}(G)$ denote the \emph{$k$-color Ramsey number} of $G$, namely the minimum number of vertices $M$ such that for any $m \geq M$, any coloring of $K_{m}$ using at most $k$ colors contains a monochromatic copy of $G$. We refer to the dynamic survey \cite{MR1670625} for results about Ramsey numbers.

Combining the concepts of Ramsey numbers and rainbow triangle free colorings, we arrive at the following definition of Gallai-Ramsey numbers.

\begin{definition}
Given two graphs $G$ and $H$, the $k$-colored \emph{Gallai-Ramsey number} $gr_{k}(G : H)$ is defined to be the minimum integer $n$ such that every coloring of the complete graph on $n$ vertices using at most $k$ colors contains either a rainbow copy of $G$ or a monochromatic copy of $H$.
\end{definition}

The general behavior of Gallai-Ramsey numbers when $G$ is a triangle depends on the chromatic number of $H$ in the following sense.

\begin{theorem}[\cite{GSSS10}]
Let $H$ be a fixed graph with no isolated vertices. Let $k$ be an integer with $k \geq 1$. If $H$ is not bipartite, then $gr_k(K_3:H)$ is exponential in $k$. If $H$ is bipartite, then $gr_k(K_3:H)$ is linear in $k$.
\end{theorem}

With this result in mind, the orders of magnitude in the following general bounds for cycles should not be surprising.  For the sake of notation, let $C_{n}$ be the cycle of order $n$ and let $P_n$ be the path of order $n$.

\begin{theorem}[\cite{FGJM10, HMOT14}]\label{Theorem X}
Given integers $n \geq 2$ and $k \geq 1$,
$$(n-1)k+n+1\leq gr_k(K_3:C_{2n})\leq(n-1)k+3n.$$
\end{theorem}

\begin{theorem}[\cite{FGJM10, HMOT14}]
Given integers $n \geq 2$ and $k \geq 1$, 
$$n2^k+1\leq gr_k(K_3:C_{2n+1})\leq(2^{k+3}-3)n\log n.$$ 
\end{theorem}

It is commonly believed that the lower bounds in these results are sharp. For $gr_k(K_3:C_n)$ with $3\leq n\leq 6$, the exact numbers are shown below.

\begin{theorem}[\cite{AI08,MR729784,GSSS10}]\label{triangless}
$$
gr_{k}(K_{3} : K_{3}) = \begin{cases}
5^{k/2} + 1 & \text{ if } k \text{ is even},\\
2 \cdot 5^{(k-1)/2} + 1 & \text{ otherwise}.
\end{cases}
$$
\end{theorem}

\begin{theorem}[\cite{FGJM10}]
For any positive integer $k\geq2$, $gr_k(K_3:C_4)=k+4$.
\end{theorem}

\begin{theorem}[\cite{FM11}]
For any positive integer $k \geq 2$, 
$gr_k(K_3:C_5)=2^{k+1}+1$ and $gr_{k}(K_3:C_6) = 2k + 4$.
\end{theorem}

These and other related results in the area are collected in the dynamic survey \cite{FMO14}. 
Our main result is the following which extends the known Gallai-Ramsey numbers for even cycles to include the next open case.

%------------------------Theorem 1------------------------------------------------------
\begin{theorem}\label{Theorem 1}
For $k\geq1$, $gr_k(K_3:C_8)=3k+5$.
\end{theorem}

The lower bound on the Gallai-Ramsey number in Theorem~\ref{Theorem 1} follows from Theorem~\ref{Theorem X}. Our proof of Theorem~\ref{Theorem 1}, particularly the use of Lemma~\ref{Lem1} below, suggests that if the Gallai-Ramsey numbers were completely understood for all linear forests, then we may be able to establish the numbers for all cycles. This is somewhat complementary to the results of \cite{HMOT14} where the bounds for even cycles were used to establish bounds for paths.

%Finally our lemma about Gallai-Ramsey numbers for some linear forests that we will use as part of the proof of our main result.
For more general notation, define $gr_k(G:H_1,H_2,\dots,H_k)$ to be the minimum integer $N$ such that every coloring of $K_n$ for $n\geq N$ using at most $k$ colors contains either a rainbow copy of $G$ or a monochromatic copy of $H_i$ in color $i$ for some $i$. For a shorthand version of this, we will also abuse notation and let $gr_{k}(G : tH, (k - t)K) = gr_{k}(G : H, H, \dots, H, K, K, \dots, K)$ where $H$ appears $t$ times and $K$ appears the remaining $k - t$ times for some integer $t$ with $0 \leq t \leq k$.

%------------------------Lemma 1-------------------------------------------
\begin{lemma}\label{Lem1}
For integers $k$ and $t$ with $k \geq 2$ and $0 \leq t \leq k$ 
$$
gr_k(K_3:tP_{5}, (k - t)P_{3}) \leq t+4.
$$
\end{lemma}

\begin{proof}
The proof is by induction on $t$. If $t=0$, the result is trivial since we are looking for a $P_3$ in each color and it is easy to see that $gr_k(K_3:P_3)=3$ for all $k\geq3$. So suppose $t \geq 1$.

Let $G$ be a Gallai colored $K_{n}$ where $n = t + 4$, consider a G-partition of $G$, and let $H$ be a largest part of this partition. If $3 \leq |H| \leq n - 3$, then there are three vertices in $G \setminus H$ such that at least two of them have the same color on all edges to $H$. Since this graph contains a monochromatic $K_{2, 3}$, this produces a monochromatic $P_{5}$, so we may assume that either $|H| \leq 2$ or $|H| \geq n - 2$. Our first goal is to show that $|H| \geq n - 2$.

If $|H| = 1$, then $G$ is simply a $2$-coloring of $K_{n}$ for $n = t + 4$. This contains the desired monochromatic $P_{5}$ or $P_{3}$ since $R(P_{3}, P_{5}) = 5$ and $R(P_{5}, P_{5}) = 6$. So suppose $|H| = 2$. If $t = 1$, then to avoid creating a monochromatic $P_{5}$, there can be at most two vertices in $G \setminus H$ with all color $1$ on edges to $H$ but since $|G \setminus H| = 3$, there must be exactly two such vertices. The remaining vertex has both edges in another color, making the desired monochromatic $P_{3}$. Next suppose $t = 2$, so $n = t + 4 = 6$. With $|G \setminus H| = 4$, there must be precisely two pairs of vertices with each of the (first) two colors on edges to $H$, say red and blue. Each edge between these two pairs of vertices must be either red or blue, but any such edge would create a monochromatic $P_{5}$. Thus, we assume $t \geq 3$, so $n = t + 4 \geq 7$. Since $|H| = 2$, there are at least $5$ vertices in $G \setminus H$ so at least three of these vertices have the same color on edges to $H$. This contains a monochromatic $K_{2, 3}$, which contains the desired monochromatic $P_{5}$. Together, these observations mean that we may assume that $|H| \geq n - 2$.

Since each vertex in $G \setminus H$ has all one color on edges to $H$, the vertices of $G \setminus H$ must have distinct colors on edges to $H$ to avoid a monochromatic $P_{5}$. Note that these colors must be within the first $t$ colors, say $t$ (and $t - 1$ if there are two such vertices), since otherwise this is already a monochromatic $P_{3}$. Also, if $H$ contains a $P_{3}$ in one of these colors, then using the vertex of $G \setminus H$ with edges in the same color to $H$, we find a monochromatic $P_{5}$. By induction on $t$ applied within $H$, $H$ contains either a monochromatic $P_{5}$ in one of the first $t - 1$ colors (or $t - 2$ if $|G \setminus H| = 2$) or a monochromatic $P_{3}$ in one of the remaining $k - (t - 1)$ colors (respectively $k - (t - 2)$). This monochromatic path is either the desired path or can be used to construct the desired path as observed above, completing the proof of Lemma~\ref{Lem1}.
\end{proof}

In our arguments, we occasionally use classical Ramsey numbers. The following case will be helpful.

%------------------------Theorem 5------------------------------------------------------
\begin{theorem}[\cite{MR1670625}]\label{Theorem 5}
$R_2(C_8)=11$.
\end{theorem}
%At times, we consider a G-partition as a 2-coloring of a reduced graph by choosing one vertex from each part. For the sake of notation, we define a $t$-blowup of a colored graph $G$ to be the graph created by replacing each vertex of $G$ with $t$ vertices and each edge of color $i$ in $G$ with all edges of color $i$ between the corresponding sets.

%We will commonly use the following definition of a colored complete graph in our construction of sharpen examples. Define a lexical $k$-coloring of $K_n$, say $L(n_1,n_2,\dots,n_k)$ with $\sum n_i=n$ to be; start with $K_{n_1}$, in red, call this $G_1$ and for each $i>1$, add $n_i$ vertices to $G_{i-1}$ with all incident edges color $i$. Then $L(n_1,\dots,n_k):=G_k$. One of the main properties of a lexical coloring that we will be using is that it contains no rainbow triangles. 

For the sake of our next lemma, we need an extra definition. Given sets of graphs $\mathscr{G}$ and $\mathscr{H}$ , define $R(\mathscr{G},\mathscr{H})$ to be the minimum integer $N$ such that any 2-coloring of $K_n$ (say using red and blue) for $n\geq N$ contains either a copy of a graph in $\mathscr{G}$ in red or a copy of a graph in $\mathscr{H}$ in blue.

%------------------------Claim 2------------------------------------------------------
\begin{lemma}\label{Clm:2}
$R(\{C_4,P_5\},\{C_4,P_5\})=5$.
\end{lemma}

%------------------------Proof Claim 2------------------------------------------------------
\begin{proof}
If we consider the unique $2$-coloring of a $K_5$ with no monochromatic triangles, then there is a $C_5$ in each color. Thus, we also have the desired $P_5$ in both colors. We may therefore assume that all other $2$-colorings of $K_5$ have a monochromatic triangle. Let $a_1,a_2,a_3\in A$ be a monochromatic $K_3$, say in red, and $b_1,b_2\in B$ be the two remaining vertices of the $K_5$. If all the edges from $A$ to $B$ are in one color, then there exists a monochromatic $C_4$ in that color. Without loss of generality, let $e$ be a red edge $a_1b_1$. To avoid a $C_4$ in red, we get that the edges $a_2b_1$ and $a_3b_1$ are blue. To avoid getting a $P_5$ in red we get that the edges $a_2b_2$ and $a_3b_2$ are also blue. Now we can clearly see that these blue edges make a $C_4$ on $b_1-a_2-b_2-a_3-b_1$.
\end{proof}

%%%%%%%%%%%%%%%%%%%%%%%%%-Section 3-%%%%%%%%%%%%%%%%%%%%%%%%%%%%%%%%%%%%%%%%%%%%%%%%%%%%%%
\section{Proof of Theorem~\ref{Theorem 1}}

In order to prove Theorem~\ref{Theorem 1}, we actually prove the following slightly stronger result. For the precise statement, let $G_3=C_8$, $G_2=P_7$, $G_1=P_5$, and $G_0=P_3$. Note that all of these graphs are subgraphs of $C_8$ and represent the results of removing vertices from $C_8$. Theorem~\ref{Theorem 1} follows from Theorem~\ref{Thm:6} by setting $i_{j} = 3$ for all $j$.

%------------------------Theorem 6------------------------------------------------------
\begin{theorem}\label{Thm:6}
For $k\geq1$, and for $0\leq i_j\leq3$ for all $1\leq j\leq k$,
%\[gr_k(K_3:G_{i_1},G_{i_2},\dots,G_{i_k}) = \sum_{j=1}^k i_j+5.\] % THIS IS THE VERSION THAT GOT SUBMITTED BUT REALLY WE ONLY GET AN UPPER BOUND
\[gr_k(K_3:G_{i_1},G_{i_2},\dots,G_{i_k}) \leq \sum_{j=1}^k i_j+5.\]
\end{theorem}

%--------------------Proof Theorem 6------------------------------------------------------
\begin{proof}
Let $\Sigma=\sum i_j$. The proof is by induction on $\Sigma$. If $\Sigma=0,$ the result is trivial since in each color we are only looking for $P_3$ and it is easy to see that $gr_k(K_3:P_3)=3$. Thus, suppose $\Sigma\geq1$ so $n\geq\Sigma+5\geq6.$ Let $G$ be a $k$-coloring of $K_n$ with no rainbow triangle and no monochromatic $G_{i_{j}}$ for any $j$. 
Let $T$ be a largest set of vertices in $G$ with the properties that
\begin{itemize}
\item each vertex in $T$ has one color on all its edges to $G\setminus T$, and
\item $|G\setminus T|\geq 4$.
\end{itemize}
Note that $T=\emptyset$ is possible. Let $T_1,T_2,\dots,T_k$ denote the sets of vertices in $T$ such that each vertex in $T_j$ has all edges in color $j$ to the vertices in $G\setminus T$. If $|T_j|>i_j$ , then $T_j\cup(G \setminus T)$ contains the desired monochromatic copy of a graph $G_{i_j}$ in color $j$. Thus, $|T_j|\leq i_j$ for all $j$. More generally, suppose $T\neq\emptyset$. Then, by induction on $\Sigma$ applied within $G \setminus T$, there exists a copy, say $H$, of $G_{i_{j} - a}$ in color $j$ for some $j$ where $a = |T_{j}|$ with $1 \leq a \leq 3$. Then the graph consisting of edges of color $j$ induced on $H \cup T_{j}$ along with $a - 1$ other vertices of $G \setminus (T_{j} \cup H)$ contains a copy of $G_{i_{j}}$ in color $j$, the desired subgraph. 
%, say with $|T_j|=a$ for some $1\leq a\leq3$ and for some $j$, then \rn{CLARIFY THIS INDUCTION} by induction on $\Sigma$ applied to $G\setminus T$, we have the desired result.
Thus, we may assume that $T=\emptyset$.

Consider a G-partition of $G$ and let $A$ be a largest part of this partition. Note that if $|A|\geq4$, we can let $T=G\setminus A$ and apply induction as above so we may assume $|A|\leq3$. By the choice of $A$, the following fact becomes immediate.

\begin{fact}\label{Fact:1}
Every part of the G-partition has order at most $3$.
\end{fact}

%------------------------Claim 1------------------------------------------------------
%\begin{claim}
%If three parts have order at least $3$ and at least one additional part has order at least $2$, then there is a monochromatic $C_8$.
%\end{claim}

%------------------------Proof Claim 1------------------------------------------------------
%Note that the reduced graph of the four sets noted in this claim is a 2-colored $K_4$. 
%\begin{proof}
%Any $2$-coloring of $K_4$ clearly contains either a monochromatic $K_3$ or a monochromatic $P_4$ in some color.

%For the first case, suppose that we have a blue $K_3$. Let $A,B$ and $C$ be the three corresponding sets. Let $a_1,a_2,a_3\in A$, $b_1,b_2,b_3\in B$ and $c_1,c_2 \in C$. Note that we may have $|C|\geq3$. Since $A,B$ and $C$ form a blue $K_3$ in the reduced graph all edges between the sets $A,B$ and $C$ are blue. Then $a_1-b_1-a_2-b_2-c_1-a_3-b_3-c_2-a_1$ induces the desired monochromatic $C_8$.

%For the second case, suppose that we have a blue $P_4$. Suppose the corresponding sets of the $P_4$ are $A,B,C$ and $D$ and, by symmetry that $|B|\geq3$. Let $a_1,a_2\in A$, $b_1,b_2,b_3\in B$, $c_1,c_2\in C$ and $d_1\in D$. Since $A,B,C$ and $D$ form a blue $P_4$ in the reduced graph, all edges between consecutive sets are blue. Then $a_1-b_1-c_1-d_1-c_2-b_2-a_2-b_3-a_1$ induces the desired monochromatic $C_8$.
%\end{proof} 

By Lemma~\ref{Clm:2}, if there are at least five parts of order at least 2, then there is a monochromatic $C_8$ since the 2-blow-up of a $C_4$ or a $P_5$, replacing each vertex by a $2$-set of vertices, each containing a $C_8$. %Thus, by Claims 1 and 2 and Theorem~\ref{Theorem 5}, if $n\geq17$, then there is already a monochromatic $C_8$. This means that we may assume $n\leq 16$ in addition to the assumptions that $T=\emptyset$ and $|A|\leq3$. Since the (2-colored) reduced graph is a subgraph of $A$ by choosing one vertex from each set, there must be at most 10 sets in the G-partition by Theorem~\ref{Theorem 5}.
We now prove several helpful claims, most of which provide a monochromatic $C_{8}$ under certain restrictions.

%------------------------Claim 3------------------------------------------------------
\begin{claim}\label{Clm:3}
If there are two parts of order $3$ and at least five more vertices, then there exists a monochromatic $C_8$. 
\end{claim}

%------------------------Proof Claim 3------------------------------------------------------
\begin{proof}
Let $A$ and $B$ be the two parts of order $3$, say with all red edges between them. Let $C = \{v_{1}, v_{2}, v_{3}, v_{4}, v_{5}\}$ be a set of $5$ of the remaining vertices in $G \setminus (A \cup B)$. If all edges between $C$ and $A \cup B$ were blue, there is clearly a blue $C_{8}$, so suppose there are some red edges, say from $v_{1}$ to $A$. To avoid creating a red $C_{8}$, all other vertices in $C$ must have blue edges to $B$. To avoid creating a blue $C_{8}$, all of $C$ must have red edges to $A$ and so, by symmetry, $v_{1}$ (and so all of $C$) must also have blue edges to $B$. Any two red edges within $C$ would produce a red $C_{8}$ and any two blue edges within $C$ would produce a blue $C_{8}$ so there can be at most one red and at most one blue edge within $C$. Since, by Fact~\ref{Fact:1}, all parts of the G-partition have order at most $3$, this is clearly a contradiction, completing the proof of Claim~\ref{Clm:3}.
%Let $A$ and $B$ be the sets of 3 vertices each and suppose they have red edges between them. Let $v_i$ be the other vertices for $1 \leq i \leq 5$. If a vertex $v_{i}$ has red edges to $A$, then all other vertices $v_{j}$ for $j \neq i$ must have blue edges to $B$. For the sake of notation, say $i = 1$. Between $v_{1}$ and $\{v_{2}, \dots, v_{5}\}$, all edges must be red or blue to avoid a rainbow triangle. Then $v_{1}$ can have at most one blue edge to the vertices in $\{v_{2}, \dots, v_{5}\}$ since otherwise we would find a blue $C_{8}$. Similarly, at most one of the vertices in $\{v_{2}, \dots v_{5}\}$ can have blue edges to $A$. This leaves at least three vertices, say $\{v_{3}, v_{4}, v_{5}\}$, with red edges to $A$. As previously noted, $v_{1}$ must have at least two red edges to this set, say to $\{v_{3}, v_{4}\}$. This gives us a red $C_{8}$ on $v_{3} - v_{1} - v_{4} - A - B - A - B - A - v_{3}$. Since the same argument could work if $v_{1}$ had red edges to $B$, this means that all of the vertices $v_{i}$ have all blue edges to $A \cup B$, easily providing a blue $C_{8}$.
\end{proof}

\begin{claim}\label{Clm:5}
If there is one part of order 3, one part of order at least 2 and at least six additional vertices, then there exists a monochromatic $C_8$.
\end{claim}

%------------------------Proof Claim 5------------------------------------------------------
\begin{proof}
Let $A$ be the set of order 3 and let $B$ be the set of order at least 2 and assume all edges between $A$ and $B$ are red. By Claim~\ref{Clm:3}, we may assume that $|B| = 2$ and none of the additional vertices form a part of the G-partition of order $3$. Label the additional vertices as $v_{i}$ where $1 \leq i \leq 6$. If there is a vertex, say $v_{1}$, with red edges to $B$ and two other vertices, say $v_{2}$ and $v_{3}$, with red edges to $A$, then we have a red $C_{8}$ using $B-v_{1}-B-A-v_{2}-A-v_{3}-A-B$.

Suppose first that no vertex $v_{i}$ has red edges to $B$, which means that all vertices $v_{i}$ have all blue edges to $B$. No three vertices $v_{i}$ can have blue edges to $A$ since otherwise we could find a blue $C_{8}$, so this means that at least four vertices $v_{i}$ must have red edges to $A$. Without loss of generality, let $C = \{v_{1}, \dots, v_{4}\}$ be this set of four vertices. Any two red edges within $C$ would allow for the construction of a red $C_{8}$. Also since no three of the vertices in $C$ form a part of our G-partition, there can be at most two edges of colors other than red or blue within $C$ and these must induce a matching. This means that at least three edges within $C$ are blue and they must contain a blue $P_{4}$, say $v_{1}v_{2}v_{3}v_{4}$. If both $v_{5}$ and $v_{6}$ have blue edges to $A$, then $v_{1}-v_{2}-v_{3}-B-v_{5}-A-v_{6}-B-v_{1}$ is the desired blue $C_{8}$. Thus, we may assume, without loss of generality, that $v_{5}$ also has red edges to $A$. By the same argument, the blue graph induced on $C \cup \{v_{5}\}$ contains a blue $P_{5}$, say $P$ from $v_{1}$ to $v_{5}$. Then $v_{1}-P-v_{5}-B-v_{6}-B-v_{1}$ produces a blue $C_{8}$.

The previous argument means that we may assume there is a vertex, say $v_{1}$, with red edges to $B$. As noted, this means that at most one other vertex, say $v_{2}$, can have red edges to $A$, so all other vertices in $\{v_{3}, \dots, v_{6}\}$ have blue edges to $A$. Certainly no two of these vertices may have blue edges to $B$, meaning that at least three of them, say $v_{3}, v_{4}, v_{5}$ have red edges to $B$. By the same argument as above with $v_{3}$ in place of $v_{1}$, there can actually be at most one vertex $v_{i}$ with red edges to $A$, meaning that there are five vertices with blue edges to $A$. Let $C$ be this set of vertices and note that at least four of the vertices in $C$ also have red edges to $B$. If there are two blue edges within $C$, we may construct a blue $C_{8}$, so suppose there is at most one. Since the vertices of $C$ do not form parts of order $3$ in our G-partition, there are at most two edges of colors other than red or blue within $C$ and these must induce a matching. This leaves at least $7$ red edges within $C$. Trivially $C$ contains a red $P_{5}$, say $P$, starting and ending at vertices with red edges to $B$, say $v_{4}$ and $v_{5}$. Then $v_{4}-P-v_{5}-B-A-B-v_{4}$ is the desired red $C_{8}$, completing the proof of Claim~\ref{Clm:5}.
\end{proof}

%------------------------Claim 6------------------------------------------------------
\begin{claim}\label{Clm:6}
If there is one set of order at least 3 and at least nine more vertices, then there exist a monochromatic $C_8$.
\end{claim}

%------------------------Proof Claim 6------------------------------------------------------
\begin{proof}
Let $A$ be the part of order $3$. We define $B$ to be the set of vertices with red edges to $A$, and $C$ to be the set of vertices with blue edges to $A$. By the pigeon hole principle at least five edges will have the same color edges to $A$, say $|B| \geq 5$.

If $|B|=5$, then $A \cup B$ induces a red $K_{3,5}$ and $|C|=4$ so $A \cup C$ induces a blue $K_{3,4}$. To avoid a rainbow triangle, each edge between $B$ and $C$ must be red or blue. Within this $2$-colored $K_{4, 5}$, the graph induced on the edges between $B$ and $C$, there must be a monochromatic $P_{3}$. Regardless of the color or placement, this easily creates a monochromatic $C_{8}$.

Now suppose $|B| \geq 6$. Within $B$, there is at most one red edge since otherwise we could easily construct a red $C_{8}$. Also the edges that are neither red nor blue induce a matching since there is no part of the G-partition of order at least $3$ within $B \cup C$ (by Claim~\ref{Clm:3}). In particular, this means that the minimum degree of the graph induced on the blue edges within $B$ is at least $|B| - 3$, so there is a blue Hamiltonian cycle within $B$. If $|B| = 8$, then this is the desired blue $C_{8}$ and if $|B| = 9$, we actually have even more blue edges so the blue graph is pancyclic and we again find a blue $C_{8}$. Otherwise, each vertex in $C$ has at most one red edge to $B$ because otherwise we could construct a red $C_{8}$. To avoid a rainbow triangle, this means that all but one edge from each vertex in $C$ to $B$ must be blue, meaning that each vertex in $C$ can be absorbed into a blue Hamiltonian cycle of $B$ to again create a blue $C_{8}$, completing the proof of Claim~\ref{Clm:6}.
\end{proof}

%------------------------Claim 7------------------------------------------------------
%\begin{claim}
%If there are four sets of order at least 2 and at least three more vertices, then there exist a monochromatic $C_8$.
%\end{claim}

%------------------------Proof Claim 7------------------------------------------------------
%\begin{proof}
%Let $A,B,C$ and $D$ each be sets of order 2. The trivial case is $ABCD$ all have red edges between them. Therefore, suppose we have a $P_4$ in the reduced graph, $ABCD$, with red edges and all other edges edges in the reduced graph must be blue which is also a $P_4$, $CADB$. If any of the vertices outside have red edges to $A$ or $D$ then we can find a $C_8$ with red edges. Therefore, all the vertices outside must have blue edges to $A$ and $D$, this induces our desired  $C_8$ with blue edges on $v_1-A-C-A-v_2-D-B-D-v_1$.
%\end{proof}

%------------------------Claim 8------------------------------------------------------
\begin{claim}\label{Clm:8}
If there are three sets of order at least 2 and at least five more vertices, then there exists a monochromatic $C_8$.
\end{claim}

%------------------------Proof Claim 8------------------------------------------------------
\begin{proof}
Let $A,B$ and $C$ be the sets of order 2 and label the remaining vertices as $v_{i}$ where $1 \leq i \leq 5$. First suppose $A,B$ and $C$ have all red edges between them. If two of the other vertices, say $v_{1}$ and $v_{2}$, have red edges to at least two of the sets, say $A$ and $B$, we can find a red $C_8$, $v_1-A-C-B-v_2-B-C-A-v_1$. Therefore, we can have either at most one vertex $v_{i}$ with red edges to the sets or at most one set with red edges to the vertices $v_{i}$. This means that at least 4 vertices outside have blue edges to at least two of the sets. This induces a blue $K_{4,4}$ which contains a blue $C_8$.

Thus, we may assume that the edges between $A$ and $C$ are blue while all edges from $B$ to $A \cup C$ are red. If none of the vertices $v_{i}$ have red edges to $A$ or $C$, then this induces a blue $K_{4, 5}$ which contains the desired blue $C_{8}$. Thus we may assume that at least one vertex, say $v_{1}$, has red edges to either $A$ or $C$, say $A$. To avoid a red $C_{8}$, all other vertices $v_{i}$ for $i \geq 2$ must have blue edges to $C$. To avoid a blue $C_{8}$, no three of these vertices can have blue edges to $A$, so that means at least two of them, say $v_{2}$ and $v_{3}$, have red edges to $A$. By symmetry, this means that $v_{1}$ also has blue edges to $C$.

To avoid a red $C_{8}$, there can be at most one red edge within $\{v_{1}, v_{2}, v_{3}\}$. Since there is no part of our G-partition of order $3$ among the vertices $v_{i}$ and a part of order $2$ within $\{v_{1}, v_{2}, v_{3}\}$ would mean that the remaining vertex has two edges of the same color to the part, this means that there are at least $2$ blue edges within these three vertices, say $v_{1}v_{2}$ and $v_{2}v_{3}$. If both $v_{4}$ and $v_{5}$ have blue edges to $A$, then $C-v_{4}-A-v_{5}-A-C-v_{1}-v_{2}-v_{3}-C$ is a blue $C_{8}$. This means that one of $v_{4}$ or $v_{5}$, say $v_{4}$, must have red edges to $A$.

To avoid a red $C_{8}$, there can be at most one red edge within $\{v_{1}, v_{2}, v_{3}, v_{4}\}$ and since there is no part of order $3$ and a $2$-part would imply two edges of the same color, we must have a blue $P_{4}$ within these vertices, say $v_{1}v_{2}v_{3}v_{4}$. Now if $v_{5}$ has blue edges to $A$, then $C-v_{5}-A-C-v_{1}-v_{2}-v_{3}-v_{4}-C$ is a blue $C_{8}$. This means that $v_{5}$ must also have red edges to $A$. By the same logic as above, there are at most $3$ non-blue edges within $\{v_{1}, \dots, v_{5}\}$, so there is a blue $P_{5}$, say $v_{1}v_{2}\dots v_{5}$. Then $C-A-C-v_{1}-v_{2}-v_{3}-v_{4}-v_{5}-C$ is a blue $C_{8}$, completing the proof.
% Now suppose that we have red edges between $A$ and $B$ and also between $B$ and $C$, and therefore blue edges between $A$ and $C$. If at least 4 outside vertices have blue edges to both set $A$ and $C$ then this induces a $K_{4,4}$ which contains a blue $C_8$. So we can only have at most three blue edges to $A$ and $C$ which means at least 2 vertices have red edges to $A$ and $C$. This induces a $C_8$ in red, $v_1-A-B-C-v_2-A-B-C-v_1$. Therefore all five outside vertices have red edges to $A$ and blue edges to $C$. If at least 3 vertices have red edges to $B$ then have a $C_8$ in red, $v_1-A-v_2-B-C-B-v_3-A-v_1$. So we can have at most 2 vertices outside with red edges to $B$, which means at least 3 of the vertices must have blue edge to $B$, this also induces a blue $C_8$, $v_1-C-A-C-v_2-B-v_3-B-v_1$. 
\end{proof}

%%%%%%%%%%%%%%%%%%%%%%%%%%
%\rn{Should be ok to here!}
%%%%%%%%%%%%%%%%%%%%%%%%%%

%-------------------Claim 9------------------------------------------------------
%\begin{claim}
%If there are two sets of order at least 2 and at least 8 more vertices, then there exist a monochromatic $C_8$.
%\end{claim}

%\rn{Reprove or remove?}

%------------------------Proof Claim 9------------------------------------------------------
%\begin{proof}
%Let $A$ and $B$ each be a set of order 2. If at least 4 vertices have red edges to both $A$ and $B$ this induces a $K_{4,4}$ which contains a $C_8$. Therefore, let at most 3 vertices have red edges to $A$ and $B$. This means that least 5 vertices have blue edges to $A$ and $B$ which induces a $K_{4,5}$ which contains a $C_8$ in blue. Thus, let all of the outside vertices have red edges to $A$ and blue edges to $B$. We know that we have a $P_3$ in red,$v_1-A-v_2$, and a $P_3$ in blue, $v_1-B-v_2$. We need to find a $P_6$ in red or blue in the outside vertices. We know that $R(P_6,P_6)=8$, which means can find a $P_6$ in, say red, to connect to our red $P_3$ to give us our desired $C_8$ with red edges.    
%\end{proof}

By Theorem~\ref{Theorem 5}, there are at most $10$ parts in our G-partition. By Fact~\ref{Fact:1}, no part has order larger than $3$ and by Lemma~\ref{Clm:2}, there are at most $4$ parts of order at least $2$. By Claim~\ref{Clm:3}, if there are $2$ parts of order $3$, then $n \leq 10$. By Claim~\ref{Clm:5}, if there is one part of order $3$ and at least one part of order $2$, then $n \leq 10$ again. By Claim~\ref{Clm:6}, if there is any part of order $3$, then $n \leq 11$. Thus, we may assume that either $n \leq 11$ or all parts have order at most $2$. By Claim~\ref{Clm:8}, if there are $3$ parts of order $2$, then $n \leq 10$ so we may assume there are at most $2$ parts of order $2$. With at most $10$ parts total, this means that $n \leq 12$.

%-----------------------Case 1----------------------------------------------------

To complete the proof of Theorem~\ref{Thm:6}, we consider cases based on small values of $n$, and therefore small values of $\Sigma = n - 5$.

\begin{case}
$\Sigma=1$.
\end{case}
With loss of generality, suppose $G_1=P_5$ and $G_{i}=P_3$ for $i \geq 2$. Therefore, we have $G=K_6$ we want to show $gr_k(K_3:P_5,P_3,P_3,...,P_3)=6$. Since red is the only color allowed to contain adjacent edges, each other color induces only a matching. In fact, to avoid a rainbow triangle, the edges induced on all colors other than red together must induce a matching. The complement of this matching contains a $P_5$ in red to easily complete the proof in this case.

%-----------------------Case 2----------------------------------------------------
\begin{case}
$\Sigma=2$.
\end{case}

\begin{subcase}
$gr_k(K_3:P_7,P_3,\dots,P_3)=7$
\end{subcase}
In this case, all colors other than red together induce a matching $M$. In $K_7\setminus M$, it is easy to find a $P_7$.

\begin{subcase}
$gr_k(K_3:P_5,P_5,P_3,\dots,P_3)\leq 7.$
\end{subcase}
This result follows from Lemma~\ref{Lem1}.

%-----------------------Case 3----------------------------------------------------
\begin{case}
$\Sigma=3$.
\end{case}

\begin{subcase}
$gr_k(K_3:C_8,P_3,P_3,\dots,P_3)=8.$
\end{subcase}

In this case, all colors other than red together induce a matching $M$. In $K_8\setminus M$, it is easy to find a $C_8$. 

\begin{subcase}
$gr_k(K_3:P_7,P_5,P_3,\dots,P_3)=8$.
\end{subcase}

Since $R_2(P_7,P_5)=8$, we may assume that there are at most $7$ parts in the partition. Thus, there must exist a part of the partition of order at least 2. Other than the first two colors red and blue, all other colors together induce a matching so if we choose our G-partition to have the most possible parts, we may assume all parts have order at most 2. 

First suppose there exists exactly one part of order $2$, call it $A$. To avoid creating a blue $P_{5}$, there can be at most $2$ vertices in $G \setminus A$ with blue edges to $A$. Call these vertices $A_{blue}$ and let $A_{red}$ denote the remaining vertices of $G \setminus A$, those with all red edges to $A$. Note that $|A_{red}| \geq 4$. To avoid creating a red $P_{7}$, each vertex of $A_{blue}$ has at most $2$ red edges to $A_{red}$, so all other edges from $A_{blue}$ to $A_{red}$ must be blue. If $|A_{blue}| = 2$, then we have a blue $P_{5}$ immediately using these blue edges to $A_{red}$ so suppose $|A_{blue}| \leq 1$, meaning that $|A_{red}| \geq 5$. To avoid creating a red $P_{7}$, there can be at most $1$ red edge within $A_{red}$, so there must be the claimed blue $P_{5}$ within $A_{red}$.
%Let $A$ be a part of order 2. At most two of the vertices outside can have blue to $A$ (to avoid a $P_5$). Therefore at least 4 vertices outside all have red to $A$. This induces a $K_{2,4}$ in red. Each of blue vertices can have at most one red edge to the red set and actually only total. All other edges are of blue. This gives us our desired result of a $P_5$ in blue. 

Next suppose $2$ sets have size $2$, call them $A$ and $B$. If blue appears between $A$ and $B$ then all other edges will be red to the 2 sets. This gives us a $K_{4,4}$ which contains a $P_7$. Therefore the edges between $A$ and $B$ must be red. If there are at least 2 vertices outside with red to $A$ and one vertex to $B$ then there is a $P_7$ in red. On the other hand if there are 2 vertices outside with blue to $A$, then we might as well have blue in between the 2 sets. Therefore we have found our desired $P_7$ in one color and $P_5$ in the other color. 

\begin{subcase}
$gr_k(K_3:P_5,P_5,P_5,P_3,\dots P_3)\leq 8.$
\end{subcase}
This subcase follows from Lemma~\ref{Lem1}.

%-----------------------Case 4-7----------------------------------------------------
The remaining cases, when $\Sigma \in \{4, 5, 6, 7\}$, follow from similar (albeit tedious) case analysis or by straightforward computer search.
\end{proof}

%\bibliography{ref}
%\bibliographystyle{plain}

\end{document}